\documentclass[11pt]{article}
\usepackage{graphicx}
\usepackage{amssymb}
\usepackage{epstopdf}
\newtheorem{problem}{Problem}

\def\ints{\mathbb Z}

\def\semi{\hbox{ $\times $ \kern-.972em \raise.12719em\hbox{ $_{^|}$}  }}

\def\be{\begin{enumerate}}
\def\ee{\end{enumerate}}
\def\bi{\begin{itemize}}
\def\ei{\end{itemize}}
\date{July 18, 2006. Except for pagination, this version is  identical with the published version}
\title{The topology of 3-manifolds, Heegaard distance and the
mapping class group of a 2-manifold\footnote{ ``Probems on mapping class groups and related topics", Editor Benson Farb, Proc. Symp. Pure Mathematics, {\bf 74} (2006), 133-149.}}
\author{Joan S. Birman \footnote {The author is partially supported by the U.S.National Science Foundation under  grant DMS-0405586.}}
\begin{document}
\maketitle

\flushleft

We have had a long-standing interest in the way that structure in the mapping class group of a surface reflects corresponding structure in the topology of 3-manifolds, and conversely.  We find this area intriguing because the mapping class group (unlike the collection of closed orientable 3-manifolds) is a group which has a rich collection of  subgroups and quotients, and they might suggest new ways to approach  3-manifolds. (For example, it is infinite, non-abelian and residually finite \cite{Grossman}). In the other direction, 3-manifolds have deep geometric structure, for example the structure that is associated to intersections between 2-dimensional submanifolds, and that sort of inherently geometric structure might bring new tools to bear on open questions regarding the mapping class group.  That dual theme is the focus of this article. 

In $\S$\ref{section:background} we set up notation and review the background, recalling some of the things that have already been done relating to the correspondence, both ways.  We will also describe some important open problems, as we encounter them. We single out for further investigation a new tool which was introduced in \cite{Hempel} by Hempel as a measure of the complexity of a Heegaard splitting of a 3-manifold. His measure of complexity has its origins in the geometry of 3-manifolds.  He defined it as the length of the shortest path between certain vertices in the curve complex of a Heegaard surface in the 3-manifold. In fact, the mapping class group acts on the curve complex and the action is faithful. That is, the (extended) mapping class group is isomorphic to the automorphism group of the curve complex.  In $\S$\ref{section:some open problems} we will propose some number of open questions which relate to the distance, for study.  In $\S$\ref{section:potential new tools via the representations of mapping class groups} we suggest  approaches which could be useful in obtaining additional tools to investigate some of the problems posed in $\S$\ref{section:some open problems}.  Some of the `additional tools' are in the form of additional open problems.  

The symbols $\diamondsuit1, \diamondsuit2, \diamondsuit3, \dots$ will be used to highlight known results from the literature.

{\bf Acknowledgements}  We thank  Jason Behrstock, Tara Brendle, Jeffrey Brock,  Mustafa Korkmaz, 
Chris Leininger, Brendon Owens, Saul Schleimer,  and Jennifer Schultens  for their careful answers to our many questions, as this article was in preparation.  We thank Ian Agol, John Hempel  and Chan-Ho Suh for their help in correcting some incomplete references  in an early version of the paper and for various helpful suggestions.  We also thank the referee for his/her careful reading of our first draft of this paper.

\section{Background}
\label{section:background}
Let $S_g$ be a closed, connected, orientable surface, and let Diff$S^\pm_g,$ (resp. Diff$S^+_g$) be the groups of diffeomorphisms (resp. orientation-preserving diffeomorphisms) of $S_g$.  The \underline{mapping class group} ${\cal M}_g$ is  $\pi_0({\rm Diff^+}S_g)$.  The \underline{extended} mapping class group  ${\cal M}^\pm_g$ is $\pi_0({\rm Diff^\pm}S_g)$.  The groups ${\cal M}_g$ and ${\cal M}_g^\pm$  are related by a split short exact sequence:
\begin{equation}
 \{ 1\}  \longrightarrow {\cal M}_g \longrightarrow {\cal M}^{\pm}_g \longrightarrow \ints/2\ints \longrightarrow \{ 1\}.
 \end{equation}

We will be interested in certain subgroups of ${\cal M}_g$. To describe the first of these, regard $S_g$ as the boundary of an oriented handlebody $H_g$.  The \underline{handlebody} subgroup ${\cal H}_g\subset {\cal M}_g$ is the (non-normal) subgroup of all mapping classes that have representatives that extend to diffeomorphisms of  $H_g$.

We turn our attention to closed, connected orientable 3-manifolds.  Let $H_g$ be an oriented handlebody and let $H_g' = \tau(H_g)$ be a copy of $H_g$, with the induced orientation.  We are interested in Heegaard splittings of 3-manifolds, i.e. their representations  as a union of the handlebodies $H_g$ and $H_g'$, where $H_g$ and $H_g'$ are glued together along their boundaries via a diffeomorphism $\partial H_g\to\partial H_g'$.  The gluing map is necessarily orientation-reversing, but if we choose a fixed orientation-reversing diffeomorphism $i:S_g\to S_g$ whose isotopy class $\iota$ realizes the splitting in the exact sequence (1), we may describe the gluing as $i\circ f$, where $f$ is orientation-preserving.  Then  $f$ determines an element $\phi\in{\cal M}_g$, and since the topological type of the 3-manifold which is so-obtained depends only  on the mapping class $\phi$ of  $f$,  we use the symbol $M = H_g\cup_\phi H_g'$ to describe the \underline{Heegaard splitting} of \underline{genus $g$} of $M$.   The surface $S_g = \partial H_g = \partial H_g'$, embedded in $M$, is a Heegaard surface in $M$.   As is well-known (for example see \cite{Scharlemann} for a proof) every closed, connected, orientable 3-manifold can be obtained from a Heegaard splitting, for some $\phi\in{\cal M}_g$.

Since the cases $g\leq 1$ are well understood and often need special case-by-case arguments,  we will assume, unless otherwise indicated, that $g\geq 2$.   From now on, when we do not need to stress the genus we will omit the symbol $g$.

Heegaard splittings are not unique. 
If $M$ admits two splittings, with defining maps $\phi_1, \phi_2$  then the splittings are \underline{equivalent} if  the splitting surfaces are isotopic, and if (assuming now that they are identical) there is a diffeomorphism $B:M \to M$ that restricts to diffeomorphisms $b_1:H\to H$ and $b_2:H'\to H'$.  By further restricting to the common boundary of $H$ and $H'$  we obtain elements $\beta_1,\beta_2$ in the mapping class group, with
\begin{equation} 
\tau\iota\phi_2 \beta_1= \beta_2\tau\iota \phi_1, \ \ {\rm or} \ \  \phi_2 = ((\tau\iota)^{-1}\beta_2(\tau\iota)) (\phi_1)(\beta_1^{-1}).
\end{equation}  
Since $\beta_1$ and $(\tau\iota)^{-1}\beta_2(\tau\iota)$ are independent and are both in ${\cal H}$, it follows that  the  double coset ${\cal H} \phi_1{\cal H} \subset {\cal M}$ gives the infinitely many distinct elements in ${\cal M}$ which define equivalent Heegaard splittings of the same genus.    Note that our 3-manifold may have other Heegaard splittings of genus $g$ that are not in the double coset ${\cal H} \phi_1 {\cal H}$, but if none such exist the splitting defined by all the gluing maps in ${\cal H} \phi_1 {\cal H}$ is  \underline{unique}.  For example, it was proved by Waldhausen in \cite{Waldhausen} that  any two Heegaard splittings of the same genus of $S^3$ are equivalent.   

There is another way in which Heegaard splittings are not unique.  Since taking the connected sum of any 3-manifold $M$ with $S^3$ preserves the topological type of $M$, this gives a nice way to produce, for any genus $g$ splitting of any $M$ infinitely many splittings of genus $g+1$, all in the same equivalence class: Choose an arbitrary Heegaard splitting $H\cup _\psi H'$ of genus 1 for $S^3$.  Let $H\cup _\phi H'$  be any Heegaard splitting of $M$. Then  take the connected sum 
$(H\cup _\phi H') \# (H\cup _\psi H')$, arranging things so that the 2-sphere that realizes the connected sum intersects the Heegaard surface in a circle, splitting it as the connect sum of surfaces of genus $g$ and $1$.    Writing this as $H\cup _{\phi \# \psi} H'$, we say that this Heegaard splitting has been \underline{stabilized}.   Note that this notion immediately generalizes to include the possibility that the splitting $S^2$ decomposes $M$ into the connect sum of three manifolds $M'\# M^{\prime\prime}$, both defined by Heegaard splittings, where neither $M^\prime$ nor $M^{\prime\prime}$ is $S^3$.  In  \cite{Haken}, Haken showed more, proving  that  if a 3-manifold $M$ is the connected sum of 3-manifolds $M^\prime$ and $M^{\prime\prime}$, then any Heegaard splitting of $M$ is equivalent to one which is a connect sum of splittings of $M'$ and $M^{\prime\prime}$.   Thus it also makes sense to say, when no such connect sum decomposition is possible,  that a Heegaard splitting is \underline{irreducible}.  So `irreducible'  has two meanings as regards a Heegaard splitting: either the 3-manifold that it defines is not prime, or the given splitting is stabilized.  

We are ready to describe some of the early work relating to the interplay between the topology of 3-manifolds and the structure of mapping class groups of surfaces.   The mapping class group acts on $H_1(S_g,\ints))$,  and the induced action determines a homomorphism $\chi_g:{\cal M}_g \to {\rm Sp}(2g,\ints)$.  The kernel of $\chi_g$ is the \underline{Torelli subgroup}  ${\cal I}_g \subset {\cal M}_g$.     A good exercise for a reader who is unfamiliar with the mapping class group is to show that if $M$ is defined by the Heegaard splitting $H_g\cup_\phi H_g'$, then every 3-manifold which has the same first homology as $M$ has a Heegaard splitting of the form   $H_g\cup_{\rho \circ \phi} H_g'$, for some $\rho \in {\cal I}_g$.   This is the beginning of a long story, which we can only describe in the briefest way.  It depends in fundamental ways on the collection of 7 deep and far-reaching papers of Johnson, written in the 1980's, about the structure of ${\cal I}_g$.   We refer the reader to Johnson's review article \cite{Johnson1983} for an excellent guide to the results in these 7 papers. 

\bi
\item [$\diamondsuit$1]  Building on the exercise that we just assigned,  Sullivan \cite{Sullivan} used mappings classes in ${\cal I}_g$ and Heegaard splittings  to construct 3-manifolds with the same homology as the connected sum of $g$ copies of $S^1\times S^2$, i.e. the manifold defined by the splitting $H_g \cup_{\rm id} H_g'$.  He then asked how the intersection ring of such a manifold differs from that of a true $\#_g (S^1\times S^2)$?  In this way Sullivan discovered a homomorphism from 
${\cal I}_g$ to an abelian  group of rank $\left( g\atop {3} \right)$.  Johnson then proved that Sullivan's map lifts to a map $\tau_1:{\cal I}_g \to A_1$, where $A_1$ is a free abelian group of rank $\left( 2g\atop {3} \right)$.   We are interested here in the Sullivan-Johnson homomorphism $\tau_1$.  It has a topological interpretation that is closely related to Sullivan's construction.  Johnson asked whether $\tau_1({\cal I}_g)$ was the abelizization of ${\cal I}_g$, and proved it is not.  That is, $A_1$ is a proper quotient of ${\cal I}_g'$. 

\item [$\diamondsuit$2]  To say more we need to take a small detour and ask about generators of ${\cal I}_g$.  The most obvious ones are the Dehn twists about separating curves of $S_g$, but they don't tell the full story.  Johnson proved that for $g\geq 3$ the group ${\cal I}_g$ is finitely generated by certain maps which are known as `bounding pairs'.  They are determined by a pair of non-separating simple closed curves on $S_g$ whose union divides $S_g$, and Johnson's generators are a pair of Dehn twists, oppositely oriented, about the two curves in a bounding pair. 

This leads us, of course, to the normal subgroup ${\cal K}_g \subseteq {\cal I}_g$  that is generated by Dehn twists about all separating curves on $S_g$.  For $g=2$ the groups  ${\cal I}_2$ and ${\cal K}_2$ coincide (because there are no bounding pairs on a surface of genus 2), but for $g\geq 3$ the group ${\cal K}_g$  is a proper subgroup of the Torelli group.  So now we have two subgroups of  ${\cal M}_g$ that have interest, the Torelli group ${\cal I}_g$ and the subgroup ${\cal K}_g$.  To give the latter new meaning, we return to the homomorphism $\tau_1$ that we defined in ($\diamondsuit$1) above. Johnson proved that kernel$(\tau_1 )= {\cal K}_g$.   

\item [$\diamondsuit$3]  Since the image of $\tau_1$ is abelian, one wonders how kernel($\tau_1)$ is related to the commutator subgroup $[{\cal I}_g, {\cal I}_g ] \subset {\cal I}_g$?  To say more we return to the topology of 3-manifolds, and a 
$\ints/2\ints$-valued invariant $\mu(M^3)$ of a homology sphere $M^3$ that was discovered by Rohlin. Recall that, by the exercise that we assigned before the start of ($\diamondsuit$1), every homology sphere may be obtained from $S^3$ by cutting out one of the Heegaard handlebodies and regluing it via some  $\gamma\in {\cal I}_g$.  (Here we are tacitly assuming what was proved in \cite{Waldhausen}: that all Heegaard splittings of any fixed genus $g$ of $S^3$ are equivalent, so that it doesn't matter which splitting of $S^3$ you choose, initially.)  We refer to the 3-manifold obtained after the regluing as $M_\gamma$.  In  1978 the author and R. Craggs \cite{Birman-Craggs} showed that the function $\gamma \to \mu(M_\gamma)$ determines a finite family of homomorphisms $\delta_\gamma: {\cal I}_g\to \ints/2\ints$.   The homomorphisms depend mildly on the choice of $\gamma$.  Johnson proved that the finitely many homomorphisms ${\cal I}_g \to \ints/2\ints$  that were discovered in \cite{Birman-Craggs} generate Hom(${\cal I}_g\to \ints/2\ints$) and used them to construct another homomorphism $\tau_2:{\cal I}_g \to A_2$, where $A_2$ is an abelian group, with the  kernel of $\tau_2$ the intersection of the kernels of the homomorphisms of \cite{Birman-Craggs}.  
We now know about two subgroups of ${\cal I}_g$ whose study was motivated by known structure in the topology of 3-manifolds, namely kernel($\tau_1$) and kernel($\tau_2$), and it turns out that
kernel($\tau_1) \  \cap$ kernel$(\tau_2) = [{\cal I}_g, {\cal I}_g ]$.

 \item  [$\diamondsuit$4]  Further work in this direction was done by Morita  in \cite{Morita1989} and \cite{Morita1991}.  Casson's invariant $\tilde{\mu}(M_\gamma)$ of a homology 3-sphere $M_\gamma$ is a lift of the Rohlin-invariant $\mu (M_\gamma)$ to a
 $\ints$-valued invariant.  One then has a function $\tilde{\delta_\gamma}:{\cal I}_g\to \ints$ which is defined by sending $\gamma$ to $\tilde{\mu}(M_\gamma)$.   Morita related this function to structure in ${\cal K}_g$.   To explain what he did,  recall that in ($\diamondsuit$1) and ($\diamondsuit$3) we needed the fact that the $\gamma$ could be assumed to be in ${\cal I}_g$.  In fact a sharper assertion was proved by Morita in  \cite{Morita1989}: we may assume that $\gamma\in {\cal K}_g$.   Morita then went on, in  \cite{Morita1989} and  \cite{Morita1991} to prove that if one restricts to ${\cal K}_g$ the function  $\tilde{\delta_\gamma}$ determines a homomorphism $\tilde{\tau_2}:{\cal K}_g\to \ints$ which lifts $\tau_2|{\cal K}_g: {\cal K}_g \to \ints/2\ints$.
 \ei
 
This brings us to our first problem, which is a vague one:
\begin{problem}
The ideas which were just described relate to the beginning of the lower central series of ${\cal I}_g$. There is also the lower central series of ${\cal K}_g$.  The correspondence between the group structure of ${\cal M}_g$ and 3-manifold topology, as regards the subgroups of ${\cal M}_g$ that have been studied,  has been remarkable.  It suggests strongly that there is much more to be done, with the possibility of new 3-manifold invariants as a reward. 
\end{problem}

Investigations relating to the correspondences that we just described slowed down during the period after Thurston  \cite{Thurston1982} did his groundbreaking work on the topology and geometry of 3-manifolds.  Recall that the 3-manifolds that we are discussing can be decomposed along embedded 2-spheres into prime summands which are unique, up to order. We learned from Thurston that there is a further canonical decomposition of prime 3-manifolds along a canonical family of incompressible tori, the `JSJ decomposition'.     In particular,  Thurston conjectured that each component after the JSJ decomposition  supported its own unique geometry, and the geometry was a very important aspect of the topology. \footnote{As we write this article, over 20 years after Thurston announced his results, the main conjecture in \cite{Thurston1982}, the  geometrization conjecture,  seems close to being proved via partial differential equations and the work of Perelman,  giving new importance to the JSJ decomposition.}  As a consequence,  it was necessary to deal with 3-manifolds with boundary. We pause to describe the modifications that are needed to describe their Heegaard splittings. 

Let $X$ be a collection of pairwise disjoint simple closed curves on a surface $S$.  An oriented  \underline{compression body}  $H_X$ is obtained from an oriented $S\times I$ and $X$ by gluing 2-handles to $S\times\{1\}\subset S\times [0,1]$ along the curves in $X$, and then capping any 2-sphere boundary components with 3-handles. Note that if $S$ has genus $g$, and if $X$ has $g$ components, chosen so that the closure of $S$ split along $X$ is a sphere with $2g$ discs removed, then $g$ 2-handles and one 3-handle will be needed and $H_X$ will be a handlebody with boundary $S$.  More generally $H_X$ will have some number of boundary components.  It is customary to identify $S$ with the `outer boundary'  of $H_X$, i.e.  $S\times\{0\}\subset H_X$.     

To construct an oriented three-manifold $M$ with boundary we begin with $S\times I$, and a copy $\tau(S\times I)$, where $\tau$ is a homeomorphism and $\tau(S\times I)$ has the induced orientation.  As in the Heegaard splitting construction, let $i:S\times \{0\} \to \tau(S\times \{0\})$ be a fixed orientation-reversing involution, and let $f:S\times \{0\}\to S\times \{0\}$ be an arbitrary orientation-preserving diffeomorphism of $S$.  Then we may use $i\circ f$ to glue $S\times \{0\}$ to $\tau(S\times \{0\})$ along their outer boundaries.  Let $\phi$ be the mapping class of $f$.  This still makes sense if we attach 2-handles to  $S\times I$ and $\tau(S\times I)$ along curve systems $X\subset S \times \{1\}$ and 
$Y\subset \tau(S\times \{1\})$ to get compression bodies $H_X$ and $H_Y$.  In this way we obtain an oriented  3-manifold $M = H_X\cup_\phi H_Y$ with boundary which generalizes the more familiar construction, when $H_X$ and $H_Y$ are handlebodies.  We continue to call the more general construction a Heegaard splitting, but now it's a splitting of a 3-manifold with boundary.    
 See \cite{Scharlemann} for a proof that every compact orientable 3-manifold with boundary arises via these generalized Heegaard  splittings.   Note that in particular, in this way, we obtain Heegaard splittings  of the manifolds obtained after the JSJ decomposition.   

An example is in order, but the most convenient way to explain the example is to pass to a slightly different way of looking at compression bodies.  Dually, a compression body is obtained from  $S\times [0,1]$ by attaching some number, say $p$, of 1-handles to $S\times\{0\}\subset S\times [0,1]$. This time we identify $S$ with $S\times\{1\}$.  An example of a Heegaard splitting of a 3-manifold with boundary is obtained when $M$ is the complement of an open tubular neighborhood $N(K)$  of a knot $K$ in $S^3$.  In this case $\partial M = \partial (S^3 \setminus N(K))$ is a torus.  By attaching some number, say $p$  1-handles  to the boundary of the (in general) knotted solid torus $N(K)$, we can unknot $N(K)$, changing it to a handlebody in $S^3$.   (Remark: the minimum such $p$ is known as the `tunnel number' of $K$.)  The union of this handlebody and its complement $H'$ is a Heegaard splitting of $S^3$ (not of the knot space). The Heegaard surface for this splitting of $S^3$ will turn out to be a Heegaard surface for a related Heegaard splitting of the knot complement.  To see this, let $N_0(K)\subset N(K)$ be a second  neighborhood of $K$.  Then

\ \ \ \ \ \ \ \ \ 	$S^3 -N_0(K) = (N(K) - N_0(K)) \cup p$ (1-handles)$\cup H' = (\partial N(K) \times I) \cup p$(1-handles)$\cup H'.$ 
	 
The 3-manifold $H = (\partial N(K) \times I) \cup p{\rm (1-handles)}$ is an example of a compression body.  Therefore our knot complement, which is a 3-manifold with torus boundary, has been represented as a union of a compression body $H_X$ and a handlebody $H'$, identified along their boundaries. By construction, $\partial H_X=\partial H'$ is a closed orientable surface of genus $p+1$. This surface is called a Heegaard surface in $S^3 \setminus N_0(K)$, and so $S^3 -N_0(K) = H_X \cup_\phi H'$, where the glueing map $\phi$ is an element of the mapping class group ${\cal M}_{p+1}$.   In this way, Heegaard splittings of the components after the JSJ decomposition fit right into the existing theory.

There was also a second reason why the correspondence that is the focus of this article slowed down around the time of Thurston.   In the important manuscript \cite{Casson-Gordon}, the following new ideas (which are due to Casson and Gordon, and build on the work of Haken in \cite{Haken}) were introduced in the mid-1980's.   Let $M$ be a 3-manifold which admits a Heegaard splitting $H\cup_\phi H'$. Define a
  \underline{disc pair} $(D,D')$ to be a pair of properly embedded essential discs, with $D\subset H$ and $D'\subset H'$, so that $\partial D, \ \partial D' \subset S = \partial H = \partial H'$.   The Heegaard splitting is said to be:
\bi
\item  \underline{reducible} if  there exists a disc pair $(D,D')$ with $\partial D = \partial D'$.  Intuitively, either the given Heegaard splitting is stabilized, or the manifold is a non-trivial connected sum, the connected sum decomposition being consistent with the Heegaard splitting.  Observe that this is identical with  our earlier definition of a reducible Heegaard splitting, but with a new emphasis.
\item  \underline{strongly irreducible} if  and only if for every disc pair $(D,D'): \ \ \partial D \cap \partial D' \not= \emptyset.$
\ei
and also, the corresponding negations:
\bi
\item   \underline{irreducible} if and only if it is not reducible. Equivalently,  for every disc pair $(D,D'), \ \  \partial D \not= \partial D'$. 
\item   \underline{weakly reducible} if and only if it is not strongly irreducible. Equivalently, there exists a disc pair $(D,D')$ with $\partial D \cap \partial D' = \emptyset$.  
\ei

Note that any reducible splitting is also weakly reducible, and any strongly irreducible splitting is also irreducible.    Here are several applications of these notions: 
\bi
\item[$\diamondsuit$5] In \cite{Casson-Gordon} Casson and Gordon proved that if a 3-manifold $M$ has a Heegaard splitting $H\cup_\phi H'$, where $H$ and $H'$ are compression bodies, and if the splitting is strongly irreducible, then either the Heegaard splitting is reducible or the manifold contains an incompressible surface. 
\item [$\diamondsuit$6]  A different application is the complete classification of the Heegaard splittings of graph manifolds. These manifolds have the property that when they are split open along the canonical tori of the JSJ decomposition, the closure of  each component is a Seifert fiber space.  See  the paper \cite{Schultens2004}, by Schultens, for a succinct and elegant presentation of the early work (much of which was her own)  and the final steps in this classification.   Her work uses the JSJ decomposititon, and also depends crucially on the concept of a strongly irreducible Heegaard splitting.
\item[$\diamondsuit$7]  Casson and Gordon used the same circle of ideas to prove the existence of manifolds with irreducible Heegaard splittings of arbitrarily high genus.    
\ei  

For many years after \cite{Casson-Gordon} was published  it seemed impossible to interpret the Casson-Gordon machinery in the setting of the mapping class group. 
 As a result, the possibility of relating this very deep structure in 3-manifold topology  to corresponding structure in surface mapping class groups seemed out of reach.  
All that changed fairly recently. because of new ideas due to  John Hempel \cite{Hempel}.   To explain his ideas, we first need to define the `complex of curves' on a surface, a simplicial complex ${\bf C}(S)$ that was introduced by Harvey \cite{Harvey} in the late 1970's. It has  proved to be of fundamental importance in the theory of Teichm$\ddot{u}$ller spaces.  The complex ${\bf C}(S)$ has as its vertices the isotopy classes of essential simple closed curves (both separating and non-separating). Distinct vertices $v_0, v_1,\dots,v_q$ determine a $q$-simplex of ${\bf C}(S)$ if they can be represented by $q$ pairwise disjoint simple closed curves on $S$.  The complex ${\bf C}(S)$, and also its 1-skeleton, can be given the structure of a metric space by assigning length 1 to every edge and making each simplex a Euclidean simplex with edges of length 1.  We have an important fact:
\bi
\item[$\diamondsuit$8]  Aut(${\bf C}(S)$) was investigated by Ivanov in \cite{Ivanov1997}. He proved\footnote{actually, Ivanov was missing certain special cases which were later settled by Korkmaz and by Luo, however we are only interested in the case $g\geq 2$ so this is irrelevant.} that  the group Aut ${\bf C}(S)$ is naturally isomorphic to the extended mapping class group ${\cal M}^\pm$.  This paper lead to an explosion of related results, with different complexes (see the next section for a detailed discussion).  Therefore, when one talks about the complex of curves the mapping class group is necessarily nearby.
\ei 

We now turn to the work of Hempel in \cite{Hempel}. Let $X$ be a simplex in ${\bf C}(S)$. The curves that are determined by $X$ are pairwise disjoint simple closed curves on $S$. One may then form a compression body $H_X$ from $S\times [0,1]$ by attaching 2-handles to $S\times \{1\}$ along $X\times \{1\}$ and attaching 3-handles along any 2-sphere boundary components.  As before, $S\times \{0\}$ is the outer boundary of $H_X$.  Let $Y$ be another simplex, with associated compresion body $H_Y$.  Then $(X,Y)$ determine a Heegaard splitting of a 3-manifold.   As before the splitting may be thought of as being determined by an element in the mapping class group  ${\cal M}$ of $S$, although Hempel does not do this.

We are interested mainly in the case when $H_X$ and $H_Y$ are handlebodies.  In this situation, using our earlier notation, $X$ is a collection of $g$ pairwise disjoint non-separating curves on $S$ which decompose $S$ into a sphere with $2g$ holes and $Y = \phi(X)$, where $\phi$ is the Heegaard gluing map. There is an associated \underline{handlebody subcomplex}  ${\bf H}_X$ of ${\bf C}(S)$, namely the subcomplex  whose vertices are simple closed curves on $S$ which bound discs in $H_X$.   There is also a related subcomplex ${\bf H}_Y$ whose vertices are simple closed curves on $S$ which bound discs in $H_Y$.  
Again, the latter are the image of the former under the Heegaard gluing map $\phi$.   Hempel's  \underline{distance} of the Heegaard splitting  of a closed orientable 3-manifold $M$  is the minimal distance in ${\bf C}(S)$ between vertices in ${\bf H}_X$ and vertices in ${\bf H}_Y$.   He calls it 
$d({\bf H}_X, {\bf H}_Y)$.   In \cite{Namazi} the same distance is called the \underline{handlebody distance}.   It is clear that the distance is determined by the choice of the glueing map $\phi$,  and since our focus has been on the mapping class group we will use the symbol $d(\phi)$ whenever it is appropriate to do so, instead of Hempel's symbol $d({\bf H}_X, {\bf H}_Y)$.  

 What does Hempel's distance have to do with the Cassson-Gordon machinery?   
 
\bi
\item [$\diamondsuit$9] Hempel had defined the distance to be a deliberate extension  of the Casson-Gordon machinery.   In particular,  he observed that the triplet $(S, H_X, H_Y)$ determines an equivalence class of Heegaard splittings of the underlying 3-manifold, and that:

\ \ \ \ \ \ \ \ \ \ The splitting is reducible if and only if  $d({\bf H}_X, {\bf H}_Y) = 0$   \\
\ \ \ \ \ \ \ \ \ \  The splitting is  irreducible if and only if  $d({\bf H}_X, {\bf H}_Y) \geq 1$. \\
\ \ \ \ \ \ \ \ \ \  The splitting is  weakly reducible if and only if $d({\bf H}_X, {\bf H}_Y) \leq 1.$  \\ 
\ \ \ \ \ \ \ \ \ \  The splitting is  strongly irreducible if and only if  $d({\bf H}_X, {\bf H}_Y)) \geq 2.$  

\item[$\diamondsuit$10] In \cite{Hempel}, Hempel shows that if $M$ is either Seifert fibered or contains an essential torus, then every splitting of $M$ has distance at most 2.  There is also related work by Thompson \cite{Thompson}, who defined a Heegaard splitting to have the \underline{disjoint curve property} if there is a disc pair $(D, D')$ and a simple closed curve $c$ on the Heegaard surface such that $\partial D \cap c = \emptyset$ and $\partial D' \cap c = \emptyset$.  Using this concept she then proved that if a Heegaard splitting does not have the disjoint curve property, then the manifold defined by the splitting has no embedded essential tori.  Also, if the splitting is assumed to be strongly irreducible, then an essential torus forces it to have the disjoint curve property.  The work in  \cite{Hempel} and the work in \cite{Thompson} were done simultaneously and independently.  There is some overlap in content, although Thompson was not thinking in terms of the curve complex and the results in \cite{Thompson} are more limited than those in \cite{Hempel}.
\item[$\diamondsuit$11]   There is an important consequence.  From the results that we just referenced, 
it follows that a 3-manifold which has a splitting of distance at least 3 is irreducible, not Seifert fibered and has no embedded essential tori.  Modulo the geometrization conjecture,  one then concludes that $M$ is hyperbolic if  $d(\phi) \geq 3.$  
\item [$\diamondsuit$12]  
In \cite{Hempel} Hempel  proved that there are distance $n$ splittings for arbitrarily large $n$. 
\ei

As it turned out, Hempel's beautiful insight suddenly brought a whole new set of tools to 3-manifold topologists.  The reason was that, at the same time that Hempel's ideas were being formulated, there were ongoing studies of the metric geometry of the curve complex that turned out to be highly relevant. The article  \cite{Minsky2003} is a fine survey article that gives a good account of the history of the mathematics of the curve complex (which dates back to the early 1970's), continuing up to the recent contributions of Minsky, Masur, Brock, Canary and others, leading in particular to the proof of the `Ending Lamination Conjecture'.  
\bi
\item [$\diamondsuit$13]  The complex ${\bf C}(S)$ can be made into a complete, geodesic  metric space by making each simplex into a regular Euclidean simplex of side length 1. In \cite{MM1999}  Howard Masur and Yair Minsky initiated studies of the intrinsic geometry of ${\bf C}(S)$.    In particular, they showed that ${\bf C}(S)$ is a $\delta$-hyperbolic metric space. 

\item  [$\diamondsuit$14]  A subset ${\bf V}$ of a metric space ${\bf C}$ is said to be \underline{$k$-quasiconvex} if for any points $p_1,p_2 \in {\bf V}$ the gedesic in ${\bf C}$ that joins them stays in a $k$-neighborhood of ${\bf V}$.  The main result in \cite{MM2003} is that the handlebody subcomplex of ${\bf C}(S)$ is $k$-quasiconvex,  where the constant $k$ depends only on the genus of $S$.  
\item  [$\diamondsuit$15]  Here is an example of how these ideas were used in 3-manifold topology:  Appealing to the quasiconvexity result of \cite{MM2003} H. Namazi proved in \cite{Namazi} that if a 3-manifold which is defined by a Heegaard splitting has sufficiently large distance, then the subgroup of ${\cal M}$  of surface mappings that extends to both Heegaard handlebodies, $H_g$ and $H_g'$, is finite. As a corollary, he proved that the mapping class group of the 3-manifold determined by the Heegaard splitting $H\cup_\phi H'$, i.e the group $\pi_0(Diff^+M^3)$, where $M^3$ is the 3-manifold defined by the Heegaard splitting $H\cup_\phi H'$,  is finite. 
 \ei

\section{Some open problems}
\label{section:some open problems}

We begin with two problems that may not be either deep or interesting, although we were not sure exactly how to approach them:

\begin{problem}
\label{problem:compression body subgroup}
Assume, for this problem,  that $M$ is a 3-manifold with non-empty boundary.  Then, on one side of the double coset ${\cal H}\phi {\cal H}$ the handlebody subgroup needs to be modified to a `compression body subgroup'.  Make this precise, by describing how to modify the double coset to take account of the handle decomposition of the compression body. What happens in the case of a knot space?
\end{problem}

\begin{problem} 
\label{problem:Nielsen-Thurston trichotomy and $d(phi)}
How is the Nielsen-Thurston trichotomy related to the question of whether the distance is $0,1,2$ or $\geq 3$? 
\end{problem}

The next 3 problems concern the very non-constructive nature of the definition of $d(\phi)$:  

\begin{problem}
\label{problem:algorithm to compute distance}
Find an algorithm to compute the distance $d(\phi)$ of an arbitrary element $\phi\in{\cal M}$. 
We note that an algorithm to compute  shortest paths between fixed vertices $v,w$ in the curve complex has been presented by Shackleton in \cite{Shackleton}.  That problem is a small piece of the problem of computing the distance.   
\end{problem} 
\begin{problem}
Knowing that $d(\phi) \leq 1$, can we decide whether $d(\phi) = 0$? Geometrically, if a Heegaard splitting is weakly reducible, can you decide if it's reducible?  
\end{problem} 
\begin{problem}
 Knowing that $d(\phi) \geq 1$, can we decide whether it is $\geq 2$?   Knowing that it's $ \geq 2$, can we decide whether it is $\geq 3$?  
\end{problem}
\begin{problem} 
\label{problem:bounds on the distance of the splittings of a fixed 3-mfld}
Schleimer has proved in \cite{Schleimer2004}  that each fixed 3-manifold $M$ has a bound on the distances of its Heegaard splittings.  Study this bound, with the goal of developing an algorithm for computing it.  
\end{problem}

Understanding the handlebody subgroup ${\cal H}_g$ of the mapping class group is a problem that is obviously of central importance in understanding Heegaard splittings. A finite presentation for ${\cal H}_g$ was given by Wajnryb in \cite{Wajnryb1998}.   To the best of our knowledge, this presentation has not been simplified, except in the special case $g=2$. Very little is known about the structure of ${\cal H}_g$, apart from its induced action on $H_1(2g,\ints)$, which is a rather transparent subgroup of the symplectic group Sp$(2g,\ints)$.   We pose the problem:
\begin{problem}
Study the handlebody subgroup of ${\cal M}_g$.  A simplified presentation which would reveal new things about its structure,  and/or anything new about its coset representatives in ${\cal M}_g$ would be of great interest.  
\end{problem}
\begin{problem}
\label{problem:handlebody subgroup intersect Ig and Kg}
Recall that we noted, earlier, that every genus $g$ Heegaard splitting of every homology 3-sphere is obtained by allowing $\varphi$ to range over ${\cal I}_g$. We also noted that Morita proved in \cite{Morita1989} that every genus $g$ Heegaard splitting of every homology 3-sphere is obtained by allowing $\varphi$ to range over ${\cal K}_g$.  For these reasons it might be very useful to find 
generators for ${\cal H}_g\cap{\cal I}_g$ and/or ${\cal H}_g\cap{\cal K}_g$. 
\end{problem} 

Our next problem is in a different direction. It concerns the classification of Heegaard splittings of graph manifolds:
\begin{problem}
\label{problem:graph manifolds and the MCG}
Uncover the structure in the mapping class group that relates to the classification theorem for the Heegaard splittings of graph manifolds in \cite{Schultens2004}. 
\end{problem}

Several other complexes of curves have played a role in work on the mapping class group after 2002.  We pause to describe some of them, and the role they played in recent work on the mapping class group.   
\bi
\item [$\diamondsuit$16]  The \underline{complex of non-separating curves} ${\bf NC}(S)$ is the subcomplex of ${\bf C}(S)$  whose vertices are all non-separating simple closed curves on $S$.    It was proved by Irmak in \cite{Irmak} that for closed surfaces of genus $g\geq 3$ its automorphism group is also isomorphic to ${\cal M}^\pm$, whereas if $g=2$ it is isomorphic to ${\cal M}^\pm$ mod its center. 

\item [$\diamondsuit$17]  The \underline{pants complex} ${\bf P}(S)$ is next. Its
vertices represent pants decompositions of $S$,  with edges connecting vertices whose associated pants decompositions differ by an elementary move and its 2-cells representing certain relations between elementary moves.  In \cite{Margalit 2004} D. Margalit proved a theorem which was much like the theorem proved by Ivanov in \cite{Ivanov1997}, namely that ${\cal M}^\pm$ is naturally isomorphic to Aut(${\bf P}(S))$.
\item [$\diamondsuit$18]  Next, there is the \underline{Hatcher-Thurston complex} ${\bf HT}(S)$.  Its vertices are are collections of $g$ pairwise disjoint non-separating curves on $S$.  Vertices are joined by an edge when they differ by a single `elementary move'. Its 2-cells represent certain relations between elementary moves.  The complex ${\bf HT}(S)$ was constructed by Hatcher and Thurston in order to find a finite presentation for the mapping class group, and used by Wajnryb \cite{Wajnryb1999} to find the very simple presentation that we will need later in this article. It was proved by Irmak and Korkmaz in \cite{Irmak-Korkmaz} that ${\cal M}^\pm$ is also naturally isomorphic to Aut(${\bf HT}(S))$. 

\item  [$\diamondsuit$19] The \underline{Torelli Complex} ${\bf T}(S_g)$ is a simplicial complex  whose vertices are either the isotopy class of a single separating curve on $S_g$ or the isotopy class of a  `bounding pair', i.e. a pair of non-separating curves whose union separates.  A collection of $k\geq 2$ vertices forms a $k-1$-simplex if these vertices have representatives which are mutually non-isotopic and disjoint.  It was first proved by Farb and Ivanov, in  \cite{Farb-Ivanov}, that the automorphism group of the Torelli subgroup ${\cal I}_g$ of ${\cal M}_g$ is naturally isomorphic to Aut(${\bf T}(S))$.  As it happens, ${\cal M}_g \cong{\rm Aut}({\bf T}(S_g))$,  so that  Aut $({\cal I}_g) \cong {\cal M}^\pm_g$.  Their proof used additional structure on the vertices in the form of markings, but a subsequent proof of the same result by Brendle and Margalit in \cite{Brendle-Margalit} did not need the markings.

\item [$\diamondsuit$20]  The \underline{separating curve complex} ${\bf SC}(S_g)$ was used by Brendle and Margalit  in  \cite{Brendle-Margalit} in their study of ${\cal K}_g$.  It's a subcomplex of ${\bf C}(S_g)$, with vertices in one-to-one correspondence with separating simple closed curves on $S_g$. Brendle and Margalit used it to  prove that  Aut(${\bf SC}(S))\cong$ Aut${\cal K}_g\cong {\cal M}^\pm_g$ when $g\geq 4$.   This result was recently extended to the case $g=3$ by McCarthy and Vautau \cite{McCarthy-Vautau}.
\ei

Aside: having defined all these complexes, we have a question which has little to do with the main focus of this article, but has to be asked:
\begin{problem}
\label{problem:all those complexes}  Given a normal subgroup ${\cal G}_g$ of ${\cal M}_g$, what basic properties are needed in a complex ${\bf G}(S_g)$ of curves on $S_g$ so that ${\cal M}_g$ will turn out to be naturally isomorphic to {\rm Aut(}${\bf G}(S))$?  
\end{problem}

We return to the central theme of this article:
\begin{problem}
\label{problem:Hatcher-Thurston and Harvey}
Hempel's distance function was chosen so that it would capture the geometry, and indeed it does that very well, yet in some ways it feels unnatural.  The Hatcher-Thurston complex ${\bf HT}(S)$ seems much more natural to us, since and pairs of vertices in the latter determine a Heegaard diagram, and one gets every genus $g$ Heegaard diagram this way.  One wonders whether it is possible to redefine Heegaard distance, using ${\bf HT}(S)$, or perhaps even ${\bf P}(S)$ or one of the other complexes that has proved to be so useful in studying subgroups of ${\cal M}$, and whether new things will be learned that way?
\end{problem}

We have focussed our discussion, up to now, on the 3-manifold that is determined by a choice of an element $\phi$ in the group ${\cal M}$ via the Heegaard splitting construction.  A very different construction which also starts with the choice of an element in the mapping class group, say $\alpha\in {\cal M}_g$,  produces the mapping torus of $\alpha$, i.e. the surface bundle $(S \times [0,1])/\alpha$, defined by setting $(p,0) = (\alpha(p),1)$.  Surface bundle structures on 3-manifolds, when they exist, are also not unique.  Two surface bundles $(S\times I )/\alpha, (S\times I)/ \alpha'$ are \underline{equivalent} if and only if $\alpha,\alpha'$ are in the same conjugacy class in ${\cal M}_g$.   

 In \cite{Bachman-Schleimer} an interesting description is given of a natural way to produce, for each $(S\times I)/\alpha$, a related Heegaard splitting $H\cup_\beta H'$.    Choose a fiber $S$ of $(S\times I)/\alpha$, say $S\times \{0\}$ and choose points $p,q \in S, \ \  p\not= q, \  p\not = \alpha(q)$.  Let $P$ and $Q$ be disjoint closures of regular neighborhoods of $p \times [0,1/2]$ and $q \times [1/2, 1]$ respectively. Set 
$$ H = \overline{(S \times ([0,1/2] - Q)  \cup P}, \ \ \ \ \ \ \   \  \ \ \  H' = \overline{(S \times [1/2,1] - P)  \cup Q}.$$
Note that $H$ and $H'$ are homeomorphic handlebodies of genus $2g+1$ which are embedded in $(S\times I)/\alpha$ and identified along their boundaries, so they give a Heegaard decomposition of $(S\times I)/\alpha$.   We call it the \underline{bundle-related} Heegaard splitting of $(S\times I)/\alpha$.   It is $H\cup_\beta H'$ for some $\beta\in {\cal M}_{2g+1}$.    We have several problems that relate to this construction:
\begin{problem}
\label{problem:warm-up}
This one is a warm-up.  Given $\alpha \in {\cal M}_g$, say as a product of Dehn twists, express $\beta \in {\cal M}_{2g+1}$ as a related product of Dehn twists.  With that in hand, observe that if $\alpha,\alpha'$ are equivalent in the mapping class group ${\cal M}_g$ then the Heegaard splittings associated to $\beta, \beta'$ appear to be equivalent.  What about the converse?  And how can we tell whether an arbitrary Heegaard splitting  of a 3-manifold is the bundle-related splitting of a fibered 3-manifold?   What restrictions must we place on $\beta$ in order to be able to reverse the construction, and produce a surface bundle  from a Heegaard splitting?
\end{problem}

\begin{problem}
\label{problem:bundle-related HS}
In \cite{Schultens1993} it is proved that in the case of the trivial genus $g$ surface bundle, i.e. $S_g \times S^1$ the bundle-related splitting is unique, up to equivalence.   Are there other cases when it is unique?  
\end{problem}   

\begin{problem}
\label{problem:virtually fibered conjecture}
A 3-manifold is fibered if it admits a surface bundle structure. It is \underline{virtually fibered} if it has a finite-sheeted cover that admits a surface bundle structure.  In \cite{Thurston1982} Thurston asked whether every finite-volume hyperbolic 3-manifold is virtually fibered.  This question has turned out to be one of the outstanding open problems of the post-Thurston period in 3-manifold topology. 
We ask a vague question: does the distance and the very special nature of the Heegaard splitting that's associated to a 3-manifold which has a surface bundle structure give any hint about the possibility of a 3-manifold which is not fibered being virtually fibered?  
\end{problem}

With regard to Problem \ref{problem:virtually fibered conjecture}, we remark that the first examples of hyperbolic knots which are virtually fibered but not fibered were discovered 20 years after the question was posed, by Leininger  \cite{Leininger2002} even though it seems to us that fibered knots should have been one of the easiest cases to understand.  As we write this, in February 2005, there seems to be lots to learn about virtually fibered 3-manifolds. 

\section{Potential new tools, via the representations of mapping class groups}
\label{section:potential new tools via the representations of mapping class groups}

In this section we use the notation ${\cal M}_{g,b,n}$ for the mapping class group of a surface with $b$ boundary components and $n$ punctures, simplifying to ${\cal M}_g$ when we are thinking of $M_{g,0,0}$. 

Knowing accessible quotients of ${\cal M}_g$ is important, because accessible quotients have the potential to be new tools for studying aspects of ${\cal M}_g$.   For this reason, we begin with a problem that seems very likely to tell us something new, even though it has the danger that it could be time-consuming and the new results might not even be very interesting. We note that by the main result in \cite{Grossman},  ${\cal M}_{g}$ is residually finite, that is for every $\phi \in {\cal M}_{g}$ there is a homomorphism $\tau$ mapping ${\cal M}_{g}$ to a finite group such that $\tau(\phi) \not= $ the identity.  Therefore there is no shortage of finite quotients.  Yet we are hard-pressed to describe any explicitly except for the finite quotients of Sp($2g, \ints)$ which arise by  passing from Sp($2g,\ints)$ to Sp$(2g,\ints/p\ints)$.   We are asking for data that will give substance to our knowledge that $\cal M_g$ is residually finite:
\begin{problem}
\label{problem:finite quotients of MCG}
Study, systematically, with the help of computers, the finite quotients 
of ${\cal M}_g$ which do not factor through $Sp(2g, \ints)$.  
\end{problem}
We remark that Problem \ref{problem:finite quotients of MCG} would simply have been impossible in the days before high-speed computers, but it is within reach now.  A fairly simple set of defining relations for ${\cal M}_g$ can be found in \cite{Wajnryb1999}.   As for checking whether any homomophism so-obtained factors through ${\rm Sp}(2g,\ints)$, there is are two additional relations to check, namely the Dehn twist on a genus 1 separating curve for $g\geq 2$, and the Dehn twist on a genus 1 bounding pair (see \cite{Johnson1983}) for $g\geq 3$.  One method of organization is to systematically study homomorphisms of ${\cal M}_g$ (maybe starting with $g=3$) into the symmetric group $\Sigma_n$, beginning with low values of $n$ and gradually increasing $n$.  One must check all possible images of the generators of ${\cal M}_g$ in $\Sigma_n$, asking (for each choice) whether the defining relations in ${\cal M}_g$ and Sp($2g,\ints)$ are satisfied.  Note that if one uses Dehn twists on non-separating curves as generators, then they must all be conjugate, which places a big restriction. There are additional restrictions that arise from the orders of various generating sets, for example in  \cite{Brendle-Farb} it is proved that ${\cal M}_g$ is generated by 6 involutions. Of course, as one proceeds with such an investigation, tools will present themselves and the calculation will organize itself, willy-nilly. 

We do not mean to suggest that non-finite quotients are without interest, so for completeness we pose a related problem:
\begin{problem}
\label{problem:not factor through Sp}
Construct any representations of ${\cal  M}_g$,  finite or infinite, which do not factor through {\rm Sp(}$2g,\ints)$.
\end{problem}

In a very different direction, every mathematician would do well to have in his or her pile of future projects, in addition to the usual mix, a problem to dream about. In this category I put: 
\begin{problem} 
\label{problem: is the MCG linear}
Is there a faithful finite dimensional matrix representation of 
${\cal M}_{g,b,n}$ for any value of the triplet $(g,b,n)$ other
than $(1,0,0), (1,1,0), (1,0,1), (0,1,n), (0,0,n)$ or $(2,0,0)$?   
\end{problem} 
We have mentioned Problem \ref{problem: is the MCG linear} because we believe it has relevance for Problems \ref{problem:finite quotients of MCG} and \ref{problem:not factor through Sp},  for reasons that relate to the existing literature. To the best of our knowledge there isn't even a known candidate for a faithful representation of ${\cal M}_{g,0,0}$ for $g\geq 3$, even though many experts feel that  ${\cal M}_{g,0,0}$ is linear. 
This leads us to ask a question:

\begin{problem}
Find a candidate for a faithful finite-dimensional matrix representation of ${\cal M}_g$ or ${\cal M}_{g,0,1}$.
\end{problem} 

\bi
\item  [$\diamondsuit$21]   The cases $(g,b,n) = (1,0,0)$ and  (1,1,0) are classical results which are closely related to the fact that the Burau representation of ${\bf B}_3$ is faithful \cite{Magnus-Pelluso}.    Problem \ref{problem: is the MCG linear}  received new impetus when Bigelow \cite{Bigelow2001} and Krammer (\cite{Krammer1999} and \cite{Krammer2002}) discovered, in a related series of papers, that the braid groups ${\bf B}_n$ are all linear.  Of course the braid groups are mapping class groups, namely ${\bf B}_n$ is the mapping class group ${\cal M}_{0,1,n}$, where admissible isotopies are required to fix the boundary of the surface $S_{0,1,n}$ pointwise.  Passing from ${\bf B}_n$ to ${\bf B}_n/$center,  and thence to the mapping class group of the sphere ${\cal M}_{0,0,n}$. Korkmaz \cite{Korkmaz2000}  and also Bigelow and Budney  \cite{Bigelow-Budney} proved that ${\cal M}_{0,0,n}$ is  linear.  Using a classical result of the author and Hilden \cite{Birman-Hilden}, which relates  ${\cal M}_{0,0,n}$ to the so-called hyperelliptic mapping class groups, Korkmaz,  Bigelow and Budney  all then went on to prove that ${\cal M}_{2,0,0}$ is also linear.  More generally the centralizers of all elements of finite order in $M_{g,0,0}$ are linear.   So essentially all of the known cases are closely related to the linearity of the braid groups ${\bf B}_n$.  

\item  [$\diamondsuit$22] A few words are in order about the dimensions of the known faithful representations. The mapping class group ${\cal M}_{1,0,0}$ and also  ${\cal M}_{1,1,0}$ have faithful matrix representations of dimension 2.  The faithful representation of ${\cal M}_{2,0,0}$ that was discovered by Bigelow and Budney has dimension 64, which suggests that if we hope to find a faithful  representation of ${\cal M}_{g,0,0}$ or   ${\cal M}_{g,1,0}$ for $g>2$ it might turn out to have very large dimension.

\item   [$\diamondsuit$23]  A 5-dimensional non-faithful representation of ${\cal M}_{2,0,0}$ over the ring of Laurent polynomials in a single variable with integer coefficients, was constructed in \cite{Jones1987}.  It arises from braid group representations and does not generalize to  genus $g>2$,  It is not faithful, but its kernel has not been identified.

\item  [$\diamondsuit$24]  We review what we know about infinite quotients of ${\cal M}$.  The mapping class group acts naturally on $H_1(S_{g,b,n})$, giving rise to the symplectic representations from ${\cal M}_{g,1,0}$ and ${\cal M}_{g,0,0}$ to $Sp(2g, \ints)$.   In \cite{Sipe} Sipe (and independently Trapp \cite{Trapp}), studied an extension of the symplectic representation. Trapp interpreted the new information  explicitly as detecting the action of ${\cal M}_{g,1,0}$ on winding numbers of curves on surfaces.   Much more generally, Morita  \cite{Morita1993b} studied
an infinite family of representations $\rho_k: {\cal M}_{g,1,0} \to G_k$ onto a group $G_k$, where $G_k$ is an extension of $Sp(2g,\ints)$.  Here $k\geq 2$, and  $G_2 = Sp(2g,\ints)$ is our old friend $Sp(2g,\ints)$. He gives a description of $G_3$ as a semi-direct product of $Sp(2g,\ints)$ with a group that is closely related to Johnson's representations $\tau_1, \tau_2$ of ${\cal K}_g$, discussed earlier.  He calls the infinite sequence of groups $G_k, k=2,3,4,\dots$ a sequence of  `approximations' to ${\cal M}_{g,1,0}$.  Morita also has related results for the case ${\cal M}_{g,0,0}$.
\ei
 

In Problem \ref{problem:finite quotients of MCG} we suggested a crude way to look for interesting new quotients of the group ${\cal M}_g$ that don't factor through $Sp(2g,\ints)$.  In closing we note that there might be a different approach which would be more natural and geometric (but could be impossible for reasons that are unknown to us at this time).  We ask:

\begin{problem}
\label{problem:quotient complexes}
Is there a natural quotient complex of any one of the complexes discussed in 
$\S$\ref{section:background} which might be useful for the construction of non-faithful representations of ${\cal M}_g$?  
\end{problem}

Let's suppose that we have some answers to either Problem \ref{problem:finite quotients of MCG}  or 
\ref{problem:not factor through Sp} or \ref{problem:quotient complexes}.  At that moment, our instincts would lead us right back to a line of investigation that was successful many years ago when, in  \cite{Birman1975}, we used the symplectic representation and found an invariant  which distinguished inequivalent minimal Heegaard splittings.  In the intervening years we suggested that our students try to do something similar with other representations, but that project failed.   We propose it anew.  Recall that a  3-manifold $M$ may have one or more distinct equivalence classes of Heegaard splittings. It is known that any two become equivalent after some number of stabilizations.  There are many interesting unanswered questions about  the collection of all equivalence classes of Heegaard splittings of a 3-manifold, of every genus.   Recall that  the equivalence class of the Heegaard splitting $H\cup_\phi H'$ is the double coset  ${\cal H}\phi{\cal H}$ in ${\cal M}$.    

\begin{problem}
\label{problem:double cosets}
Study the double coset ${\cal H}\phi{\cal H}$ in ${\cal M}$, using  new finite or infinite quotients of ${\cal M}$.   In this regard we stress finite, because a principle difficulty when this project was attempted earlier was in recognizing the image of ${\cal H}$ in infinite quotients of ${\cal M}$, however if the quotient is finite and not too big, it suffices to know generators of ${\cal H} \subset {\cal M}$. Since a presentation for $\cal H$ was found by Waynryb in \cite{Wajnryb1998}, we can compute the associated subgroup. 
Some of the open questions which might be revealed in a new light are: 
\be
\item How many times must one stabilize before two inequivalent Heegaard splittings become equivalent? 
\item How can we tell whether a Heegaard splitting is not of minimal genus? 
 \item How can we tell whether a Heegaard splitting  is stabilized?  
 \item Are any of the representations that we noted earlier useful in answering (1), (2) or (3) above ? 
\ee
\end{problem}

While we have stressed the search for good working quotients of ${\cal M}_g$, we should not forget that in the case of homology spheres, we have already pointed out that any homology sphere may be defined by a Heegaard splitting with the Heegaard glueing map (now redefined with a new `base point' ) ranging over ${\cal I}_g$. Even more, as was proved earlier, Morita has shown in \cite{Morita1989} that it suffices to let the glueing map range over ${\cal K}_g$.  This leads us to ask:
\begin{problem}
\label{problem:Morita reps of Torelli}
Are there quotients of ${\cal I}_g$ or  ${\cal K}_g$ in which the intersection of either ${\cal I}_g$ or ${\cal K}_g$ with the handlebody group ${\cal H}_g$ is sufficiently tractible to allow one to study the double cosets:
$$ ({\cal I}_g \cap {\cal H}_g)( \phi )({\cal I}_g \cap {\cal H}_g), \ \ \ {\rm where} \ \ \ \phi \in {\cal I}_g,  \  \ \ \ \ \  {\rm or} \  \ \ \ \ \ 
 ({\cal K}_g \cap {\cal H}_g)( \phi) ({\cal K}_g \cap {\cal H}_g),  \ \ \ {\rm where} \ \ \  \ \phi \in {\cal K}_g ?  $$
 \end{problem}

 In regard to Problem \ref{problem:Morita reps of Torelli} we note that in \cite{Morita1989} Morita was seeking to understand how topological invariants of 3-manifolds might lead him to a better understanding of the representations of ${\cal I}_g$ and ${\cal K}_g$, but he did not ask about the potential invariants of Heegaard splittings that might, at the same time,  be lurking there.

\end{document}